\documentclass[ 11pt] {article}
\usepackage{a4,  psfrag, epsfig, epsf, amscd, amsfonts, amsmath, amssymb, amstext, amsthm, enumerate}
\begin{document}

\newtheorem{satz}{Theorem}[section]
\newtheorem{defin}[satz]{Definition}
\newtheorem{bem}[satz]{Remark}
\newtheorem{hl}[satz]{Proposition}
\newtheorem{lem}[satz]{Lemma}
\newtheorem{ko}[satz]{Corollary}

\title{Equicontinuous Geodesic Flows }
\author{Christian Pries  \\ Fakult$\ddot{a}$t f$\ddot{u}$r Mathematik\\ Ruhr-Universit$\ddot{a}$t Bochum\\ Universit$\ddot{a}$tsstr. 150 \\ 44780 
Bochum \\ Germany\\ Christian.Pries@rub.de }
\maketitle

\begin{abstract}
This article is about the interplay between  topological dynamics and differential geometry. 
One could ask how many informations of the geometry are carried in the dynamic of the geodesic flow. M. Paternain proved in [6] that an expansive geodesic flow on a surface implies
that there are no conjugate points. Instead of regarding notions that describe chaotic behavior (for example expansiveness) we regard a notion that describes the stability of orbits in dynamical systems, namely equicontinuity and distality. In this paper we give a new sufficient and necessary condition for a compact Riemannian
surface to have all geodesics closed (P-manifold):\\ $(M,g)$ is a
P-manifold iff the geodesic flow 
$SM\times \mathbb{R} \to SM$ is equicontinuous. We also prove  a weaker theorem
for flows on  manifolds of dimension $3$. At the end we discuss  some properties of 
equicontinuous geodesic flows on noncompact surfaces and higher dimensional manifolds.
\end{abstract}

\section[]{Introduction}
In the complete paper all geodesics are parametrized by arc length and the geodesic flow is complete. The manifolds and the  Riemannian metrics  are $C^{\infty}$.
$\pi: TM\to M$ denotes the canonical projection. At first we summarize some facts about recurrent maps and we state  theorem [2.7] (due to Boris Kolev
and Marie-Christine P$\acute{e}$rou$\grave{e}$me) about the set of fixed points of recurrent maps on surfaces. In section 3 we study the geodesic return map. In section 4 we prove that equicontinuous geodesic flows are periodic. In section 5 we prove that if a flow without singularities on a 3-dimensional manifold admits a global Poincar$\acute{e}$ section and has enough periodic orbits, then the flow is pointwise periodic. In the last section we show that the existence of an  equicontinuous geodesic flow  on a compact manifold $M$ implies that the fundamental group is finite. Moreover, we discuss  equicontinuous geodesic flows on noncompact manifolds.

\section[]{Recurrent Behavior}

\begin{defin}

A dynamical system $(X,T)$ is called distal if \\ $\inf\{ d(xt, yt)  |  t\in T \}= 0$
implies $x=y$.\\ A system $(X,T)$ is called  equicontinuous  (regular) if for all $\epsilon > 0$ there exists an $\delta (\epsilon) > 0$ such that for all $x,y$ with 
$d(x,y) < \delta (\epsilon) $ we have $d(xt,yt)< \epsilon$ for all $t \in T$.
\end{defin}

Here is  a well-known fact about equicontinuous systems on compact metric  spaces:
\begin{satz}
An equicontinuous flow $\Phi : X\times \mathbb{R} \to X$ on a compact metric  space  is uniformly almost periodic, i.e. 
for every $\epsilon > 0$ there exists an $\tau >0 $ such that in every intervall $I$ of length $\tau$ there exists an $t\in I$ such that for all $x $ 
we have $d(\Phi (t,x), x) < \epsilon $. 
\end{satz}
Proof: See theorem 2.2 in [1].

A weaker form of uniform almost periodicity for maps is the following:

\begin{defin}
A continuous map $f$ on a  metric space $(X,d)$ is recurrent if there exists a sequence $n_k\to \infty $ such that for $\sup_{x\in X} d(f^{n_k}(x),x) \to 0$ $k\to \infty$.
\end{defin}

\begin{defin}
A  continuous map $f:X\to X$ on a metric space $(X,d)$ is called paracompact-recurrent on $Y\subset X$ if there exists a sequence $n_k\to \infty $ such that for $k\to \infty$
we have 
 $\sup_{x\in C} d(f^{n_k}(x),x) \to 0$, where $C\subset Y$ is a compact subset of $X$.
\end{defin}
Note that in the definition of paracompact-recurrence the sequence $n_k$ is fixed and does not depend on $C\subset Y$ and
 that paracompact-recurrence and recurrence are independent from the metric which defines the topology if the space $X$ is compact.

\begin{lem}
If $f$ is recurrent, then $f^m$ is recurrent.
\end{lem}

Set $s_k:= \sup_{x\in X} d(f^{n_k}(x),x)$. We conclude
$$ d(f^{n_km}(x),x) \le \sum_{i=0}^{m-1} d(f^{(m-i)n_k}(x), f^{(m-i-1)n_k}(x))$$
$$ \le  \sum_{i=0}^{m-1} 
d(f^{n_k}(f^{(m-i-1)n_k}(x)), f^{(m-i-1)n_k}(x)) \le  ms_k.$$
{\hfill $\Box$ \vspace{2mm}}

\begin{lem}
Given a  continuous map on a compact surface $S$. Let $F$ be a finite non empty subset of $Fix(f)$ (the fixed point set). If $f$ is paracompact-recurrent on $S-F$, then $f$ is 
recurrent.
\end{lem}
W.l.o.g we can suppose that our metric is induced from a Riemannian metric. 
Let $\{x_1, \cdots ,x_m \}$ denote $F$ and choose a tiny $\epsilon > 0$ and
define $$C:= S-(\bigcup_i B(x_i,\epsilon) ).$$  $c_i= \partial B(x_i,\epsilon)$
defines a simple closed curve. Choose a $K:=K(\epsilon )$ such that for all
$k> K$ we have:
$$ \sup_{x \in C} d(f^{n_k}(x),x) < 4\epsilon 
\quad\mbox{and}\quad f^{n_k}(c_i) \subset B(c_i, \frac{\epsilon}{2}).$$ 
This is possible,
since there is a compact subset in $S-F$ that contains $c_i$ and $C$.  \\ Hence $f^{n_k}(B(x_i, \epsilon))
\subset B(x_i, 2 \epsilon)$ thus $d(f^{n_k}(x),x) <  4\epsilon $ for all $x\in C^c$ and $k > K$. Consequently $d(f^{n_k}(x),x) <  4\epsilon $ for all $x \in X$ and $k > K$, i.e. $f$ is recurrent. {\hfill $\Box$ \vspace{2mm}}

The following  nice theorem can be found in [4] as theorem 1.1.
\begin{satz}
A non trivial orientation-preserving and recurrent homeomorphism of the sphere $S^2$ has exactly two fixed points.
\end{satz}

\section[]{The Geodesic Return Map }

A well-known technic to study geodesic flows on a surface is to study the geodesic return map.
Let $A$ denote the open annulus. \\ Suppose that we have an equicontinuous geodesic flow $\Phi$ on the unit tangent bundle $ SM $ of an orientable  Riemannian surface $M$. Let $\gamma$ be a simple closed geodesic. Let denote $W=SM|\gamma-T\gamma$ (the set of all
 unit tangent vectors based on $\gamma$ which are not element of $T\gamma$).
Note that  $W$ is homeomorphic to the union of two  open annulus $A_0$ and $A_1$.
We can identify $v \in A_0$ with $(x,\theta )$ where $\pi(v)=x$ and 
$\theta \in (0,1)$
 is the angle of $v$ and $\dot{\gamma}(t)$ divided by $\pi$.  \\ Since every orbit is recurrent, we can define  for our flow $\Phi$ the map  $F:A_0\to A_0 $  by 
$$ F(x,\theta )=(x_0, \theta_0) , $$
where $x_0=\pi (\Phi_{t_0}(x,\theta))$   is the next intersection point of $\{ \pi (\Phi_t(x,\theta))|t>0 \} $ with 
$\gamma$   such that $\Phi_{t_0}(x,\theta)=(x_0, \theta_0)\in A_0$.
We just simple write $F:A\to A$. $F$ can be extended to an homeomorphism of $S^2$ by two-point-compactification. In this case $F$ has two
fixed points $\{\infty\}$ and $\{- \infty\}$ as we will see.
If for $v=(x,\theta)$ we have $\theta$ near zero, the geodesic $\Phi_t(v)$
stays near $\dot{ \gamma }$ (by equicontinuity), hence $F^n(x,\theta)$ is near zero for all $n$, thus  $\{\infty\}$ and $\{- \infty\}$ are fixed points.
In this paper, we call  the  extension of $F$ the geodesic return map and denote
it by $F:S^2\to S^2$.

\begin{hl}
If  $\Phi$ is equicontinuous then $F$ is recurrent on $S^2$.
\end{hl}
 We only need to show that  $F$ is paracompact-recurrent on $S^2-(\{\infty\} \cup \{- \infty\}) $, since then we can apply lemma [2.6].

For $0< \theta_0 < \theta_1 <1 $ we set 
$$K(\theta_0 ,\theta_1)=\{ (x,v) \in A  |   \theta_0 \le v \le \theta_1 \}.  $$

\begin{lem}
For  $K:=K(\theta_0 ,\theta_1)$ and large $M:=M(\theta_0 ,\theta_1)$, the constant
$$s(\theta_0 ,\theta_1, M):=\inf\{t_1-t_0  |  -M< t_0< t_1 < M,  \Phi_{t_0}(v)\in A, 
\Phi_{t_1}(v)\in A, v\in K \}$$ is strict positive. Every geodesic 
$\gamma_{v}$ intersects $\gamma$  at least 2 times in the forward direction on the
intervall $[0,M]$ and at least 2 times in the backward direction on the
intervall $[-M,0]$ if $v\in K$.

\end{lem}
Proof: Use the compactness of $K$.{\hfill $\Box$ \vspace{2mm}}

\begin{lem}
There exists a constant $\nu(\theta_0 ,\theta_1, M, r) > 0 $ for \\ every 
$ \frac{s(\theta_0 ,\theta_1 , M)}{2} > r>0$, large $M$ and 
$K(\theta_0 ,\theta_1)$ such that the following holds: \\ If for some $w\in SM$ we have 
an $v\in K(\theta_0 ,\theta_1 )$ with $d(v,w)< \nu(\theta_0 ,\theta_1,r, M) $, then there is 
an unique $t_w$ such that $|t_w|<r$ and $\Phi_{t_w}(w) \in A$.
\end{lem}
Proof: Use the compactness of $K$.{\hfill $\Box$ \vspace{2mm}}

Given $v\in A$ then for every integer  $n$ we define $t(n,v)$ as the 
following:  \\ $t(n,v)$ is the unique element of $\mathbb{R}$ such that 
$\Phi_{t(n,v)}(v)=F^n(v)$ and 
$t(0,v)=0$. Moreover given $v\in A$ and $t\ge 0$ we define:
$$P(v,t):= \max\{n \ge 0  | t(n,v) \le t \}. $$ 
Note that we can find a neighbourhood $U_0$ of $\gamma$ looking like a  strip
(M is orientable) , therefore
 $U_0 -\gamma$ is the distinct union of two open strips $S^*$ and $S_*$.
We can define what lies in $U_0 -\gamma$ above and below $\gamma$. The area where the points of $A$ are pointing in-ward is defined to be 
$S^*$ (" above "). The other points of $U-\gamma$ are lying in $S_*$ ( " below ").

 For $j\in \{0,1\}$ choose sequences $(\theta_{j,i})_i$ with 
$$\theta_{0,i+1} < \theta_{0,i} < \theta_{1,i} < \theta_{1,i+1}$$
and converging strictly to $j$.
Set $K_i=K(\theta_{0,i} ,\theta_{1,i})$ and fix a $v_0\in \cap_i K_i$.

\begin{lem}
There exist  sequences $T_i\to \infty$, $0< \zeta_i\to 0$ and $0 < \beta_i \to 0$
such that:
\begin{enumerate}
\item $\pi(\Phi_{T_i+ \zeta_i}(v)) \in S^*$ for all $v\in K_i$
\item $\pi(\Phi_{T_i- \zeta_i}(v)) \in S_*$ for all $v\in K_i$
\item $\Phi_{T_i}(v_0) \in A$ 
\item $d(\Phi_{T_i+s}(v),v)< \beta_i$ for all $v\in K_i$ and $|s|< 2\zeta_i$
\item For every $v\in K_i$ there exists an unique intersection point of $\gamma$ and $\pi(\Phi_{T_i+s}(v))$, where $s$ ranges over $|s|< 4\zeta_i$
\end{enumerate}
\end{lem}
Proof: Choose sequences $0 < \delta_i \to 0$, $0 < s_i \to 0$ and $\alpha_i \to 0$
such that if $v\in K_i$ and $w\in SM$ with $d(v,w)< \delta_i$ then 
\begin{enumerate}[i)]
\item $\pi(\Phi_{s_i}(w)) \in S^*$ 
\item $\pi(\Phi_{-s_i}(w)) \in S_*$ 
\item $d(\Phi_{t}(w),v)<\alpha_i $ for all $|t|\le 8s_i$
\item $\pi(\Phi_{t}(w))$  has an unique intersection point with $\gamma$ where $t$ ranges over $|t|\le  8s_i$.
\item $\pi(\Phi_{[-20s_i, 20s_i]}(w)) \subset U_0$

\end{enumerate}
This can be done easy. First choose  an increasing sequence $M_i$ such that we can
define $s(\theta_{0,i+1} ,\theta_{1,i+1}, M_i)$ (and therefore $s(\theta_{0,i+1} ,\theta_{1,i+1}, M_{i+1})$ is also defined ) and a decreasing sequence $s_i \to 0$ such that 
$$32s_i< \min\{\frac{s(\theta_{0,i+1} ,\theta_{1,i+1}, M_i)}{2} , \frac{s(\theta_{0,i} ,\theta_{1,i} , M_i)}{2} \} $$ and then a decreasing sequence $\delta_i > 0$
$$\delta_i < \nu(\theta_{0,i} ,\theta_{1,i}, M_i, \frac{s_i}{2} ).$$ Moreover $\delta_i$
is chosen so small, that $d(w, K_i) < \delta_i $ implies that $\Phi_{t_w}(w)$
lies in $K_{i+1}$ ($t_w$ is from definition [3.3]).
Since $s_i, \delta_i \to 0$ and our flow $\Phi$ is defined on a compact space, surely $\alpha_i$
exists and  by regarding only large $i$ we have 
$\pi(\Phi_{[-20s_i, 20s_i]}(w)) \subset U_0$.  If $d(v,w)< \delta_i$ and $v\in K_i$ then i)  ii) and v) holds. We check that  iv) holds: \\ Note that  for $d(v,w)< \delta_i$ and 
$v\in K_i$ we have  $\Phi_{t_w}( w)\in K_{i+1}$ and 
$$[-8s_i, 8s_i] \subset [-8s_i-t_w, 8s_i +t_w] \subset [-9s_i, 9s_i],$$  
 and therefore  $\Phi_{[-8s_i, 8s_i] }(w)$ lies in $\Phi_{[-9s_i, 9s_i] }
(\Phi_{t_w}( w))$ . \\ Hence we conclude that iv) holds
since $18s_i< s(\theta_{0,i+1} ,\theta_{1,i+1} , M_i)$. \\ Now
choose a sequence $S_i\to \infty$ such that $d(\Phi_{S_i}(x),x)< \delta_i$ for all $x$
(apply theorem [2.2]).
There exists an unique $|r_i|< s_i$ such that $\Phi_{S_i+r_i}(v_0)\in A$.
Set $T_i:=S_i+r_i$ and $\zeta_i:= 2s_i$. Given $v\in K_i$ note that 
$$T_i+ \zeta_i = S_i+r_i+2s_i> S_i+s_i$$ and $|T_i+ \zeta_i -S_i|<3s_i $. From this ,  
i), iv) and v) we conclude 
 $\Phi_{S_i+s_i}(v)  \in S^*$ and therefore $\Phi_{T_i+ \zeta_i}(v) \in S^*$.
Analogously we conclude from ii), iv) and v) that $\Phi_{T_i- \zeta_i}(v) \in S_*$.
Thus 1,2 and 3 is proven. \\ If $|t|\le 2\zeta_i= 4s_i$, then $|T_i+t-S_i|< 5s_i$ thus the orbit
segment \\ $\Phi_{[T_i-2\zeta_i, T_i+2\zeta_i]}(v)$ lies in the orbit segment 
 $\Phi_{ [S_i-8s_i, S_i+8s_i]}(v)$, hence  by iv) has an unique intersection point with $\gamma$. This proves 5.
Since $2\zeta_i\to 0$ and uniformly $\Phi_{T_i}(x)\to x $ and $M$ is compact, there exists
a sequence $\beta_i$.
{\hfill $\Box$ \vspace{2mm}}

Proof of the proposition: Set $p_i:=P(v_0, T_i+\zeta_i)$. We show that for this sequence
$F$ is paracompact-recurrent. Note that  
$ P(v_0, T_i+\zeta_i)=P(v_0, T_i)$ by lemma [3.4] and $K_{i}\supset K_{j}$ if $j\le i$. Given any compact set $C\subset A$ choose $I_0$ such that $C\subset K_{I_0}$. The functions
$$ G_i: A\to \mathbb{N}$$ defined by $G_i(v)=P(v, T_i+\zeta_i)$ are constant on $K_{I_0}$ 
if $i\ge I_0$. \\ Indeed, by construction the functions $G_i$
are locally constant on $K_{i}$, hence constant $p_i=P(v_0, T_i+\zeta_i)$ on the connected set $K_i$.
Therefore $G_i=p_i$ on $K_{I_0}$ if $i\ge I_0$. \\ We have $|T_i-t(v,G_i(v))|< 2\zeta_i$ for $v\in K_i$ by construction, since  we know from lemma [3.4, 1 and 2] that  $P(v, T_i-\zeta_i )= p_i-1$. Therefore we conclude from lemma [3.4, 4)] that for all $v\in K_{I_0}$ and $i\ge I_0$ we have
 $d(F^{p_i}(v),v)=d(\Phi_{t(v,p_i)}(v), v)<\beta_i$.{\hfill $\Box$ \vspace{2mm}}

\section[]{Equicontinuous Geodesic Flows On Surfaces }

\begin{defin}
$(M,g)$ is called a P-manifold if all geodesics are closed. 
\end{defin}

The following lemma is easy to prove:
\begin{lem}
If $(M,g)$ is a P-manifold then the geodesic flow $(SM,\Phi ) $ is equicontinuous.
\end{lem}
Proof: It is a well-known fact that if $M$ is a P-manifold then the flow is periodic (compare [8]) and $M$ is compact. Let $L$ denote the smallest period. Choose for $\epsilon > 0$
an $\delta > 0$ such that $d(v,w)< \delta $ implies 
$$d(\Phi_t(v),\Phi_t(w) ) < \epsilon $$
for $|t|<2L$, hence it holds for all $t$. {\hfill $\Box$ \vspace{2mm}}

To prove our first theorem we need the following theorem.
\begin{satz}[Ballmann]
Every compact Riemannian manifold $(M, g)$ of dimension $2$ has at least three simple closed geodesics.
\end{satz}
Proof: See [2]. {\hfill $\Box$ \vspace{2mm}}

\begin{satz}
Given a compact Riemannian manifold $(M,g)$ of dimension $2$  then the following conditions are equivalent:
\begin{enumerate}
\item $M$ is a P-manifold.
\item $(SM,\Phi ) $ is equicontinuous.

\end{enumerate}
\end{satz}
 1) implies 2) by lemma [4.2]. We show that 2) implies 1). Take   a simple closed geodesic $\gamma$. 

\begin{lem}
If $Z$ is a compact set in $M-\gamma$, then there are 
$0< \theta_{0,Z} < \theta_{1,Z}< 1$ such that for every geodesic $\alpha$ intersecting 
$Z$ we have that $\Phi (Z \times \mathbb{R} ) \cap A $ is a subset of 
 $K( \theta_{0,Z} , \theta_{1,Z} )$.
\end{lem}
Proof: Choose an open set $V$ around $\gamma$ such that $V\cap Z= \emptyset$.
 Choose $\delta > 0$ such that $d(w,T\gamma )< \delta$ implies
$\gamma_w \subset V$. If  $\theta_{0,Z} $ or $\theta_{1,Z} $ would not exist, we could
find  a geodesic $\gamma_w \subset V$, but starting in $Z$.
{\hfill $\Box$ \vspace{2mm}}

\begin{lem}
Every geodesic intersects $\gamma$.
\end{lem}
Given a point $x\in M-\gamma$ and $w\in S_xM$. Choose a path-connected compact set 
$C$ in $SM$ such that 
$S_xM \subset C$ and $\pi(C)\cap \gamma= \emptyset$. Apply lemma [4.5] to $Z=\pi(C)$.
Choose a curve $\beta $ in $C$ from $w\in S_xM$ to $q\in S_xM$ where $\gamma_q$ intersects $\gamma$.
Cover $\beta $ with finitely many balls $B_i$ (say n) such that for any two vectors 
$p,u \in B_i$ we have for a large $M$ and all $t$
$$d(\Phi_t(u), \Phi_t(p))< \nu(\theta_{0,Z} - \epsilon  ,\theta_{1,Z} + \epsilon , M , b), $$
where $b= \frac{s(\theta_{0,Z} - \epsilon ,\theta_{1,Z} + \epsilon, M)}{4}$ and $ \epsilon$
is small.
By induction we conclude that if 
$\gamma_q(T_0)\in \gamma$, then $\gamma_w([T_0-(n+1)b,T_0+(n+1)b])$ 
intersects $\gamma$.
{\hfill $\Box$ \vspace{2mm}}

Proof of the theorem:
Since the lifted geodesic flow on the orientable double cover $N$ is equicontinuous, it suffices to show that the theorem holds for orientable surfaces,
otherwise we regard the orientable double cover $N$ and conclude that 
$N$ is a P-manifold. Apply the construction of $F$ in section 3 to
$(M,g)$ and $\gamma$. By theorem [4.3] and lemma [4.6] we have a 
periodic point in $y\in A$, thus for some $m$ it follows that $F^{2m}$ is an orientation-preserving homeomorphism
with three fixed points ($\{\infty\}$,$\{- \infty\}$, $y$). From   proposition [3.1] we know that $F^{2m}$ is recurrent, hence trival by theorem [2.7].
Lemma [4.6] implies that every geodesic is closed.{\hfill $\Box$ \vspace{2mm}}

Note that distality of the geodesic flow does not imply that the manifold is
a P-manifold. Take the torus $T^n$ with the standard flat metric. The flow is distal.
Indeed, for an unit vector $v\in \mathbb{R}^n$ define the vector field 
$X_v(x)=v$. Note that the solutions of these vectorfields are our lifted geodesices. \\If 
$\inf \{ d(xt, yt)  |  t\in T \}= 0$ for $x,y \in ST^n$ then the lifted geodesics $xt,yt$ are solutions of the same vectorfield, but the projected flows on $T^n$ of these vectorfields are equicontinuous, since equicontinuity implies distality (see [1]) we know that $y=x$. 
Note also that on surfaces of higher genus the geodesic flow has positive entropy,
but distal flows on compact metric spaces have always zero entropy (see [5]), thus
the only compact surfaces $M$ that can admit distal geodesic flows are $S^2$, $\mathbb{R}P^2$, $T^2$ and the Klein bottle and they do.

\section[]{ Flows On Manifolds Of Dimension 3}

\begin{defin}
A global surface of section $\Sigma$ for a $C^{\infty}$ flow $\Phi$ without singularities on a three dimensional manifold is a compact submanifold  with the following properties:
\begin{enumerate}
\item If $\Sigma$ has a boundary then its boundary components are periodic orbits.
\item The interior of the surface (stated with $\stackrel{\circ }{ \Sigma }$) is transversal to $\Phi$.
\item The orbit through a point not lying on the boundary of $\Sigma$ hits the interior
in forward and backward time.
\item Every orbit intersects $\Sigma$.
\end{enumerate}

\end{defin}

 There is a natural compactification of $\stackrel{\circ }{ \Sigma }$ to a closed surface
by collapsing the boundary components to a point. We call this unique compactification
the compactification of $\stackrel{\circ }{ \Sigma }$ to a closed surface.
If the flow is equicontinuous we can do more: \\ The return map $F: \stackrel{\circ }{ \Sigma } \to \stackrel{\circ }{ \Sigma }$
can be extended to the compactification of $\stackrel{\circ }{ \Sigma }$ to a closed surface
by defining the collapsed boundary components to be fixed points. This is well defined and can be proven as in the beginning of section 3. We call this map in this section the extended poincar$\acute{e}$ section map.

\begin{satz}
Let $\Phi$ be a $C^{\infty}$ equicontinuous flow without singularities on a three dimensional manifold that admits a global surface of section $\Sigma $. Let $X$ be the compactification of 
$\stackrel{\circ }{ \Sigma }$ to a closed surface, then the following holds:
\begin{enumerate}
\item If $X$ is  homeomorphic to the torus or the  Klein bottle
 and the extended poincar$\acute{e}$ section map on $X$ has at least one periodic point, then the flow is pointwise periodic.
\item If  $X$ is homeomorphic to the sphere 
and the extended poincar$\acute{e}$ section map on $X$ has at least three periodic points, then the flow
is pointwise periodic.
\item If $X$  is homeomorphic to the projective plane  
and the extended poincar$\acute{e}$ section map on $X$ has at least two periodic points, then the flow
is pointwise periodic.
\item If  $X$  is homeomorphic to a surface
of negative Euler characteristic, then  the flow is pointwise periodic.
\end{enumerate}
\end{satz}

Proof: The proof is quite analogue to theorem [4.4], and so
we only repeat some ideas. Indeed, the proof of theorem [4.4] was actually of topological nature. The following theorem is important:

\begin{satz}
A recurrent homeomorphism of a compact surface with negative Euler characteristic is periodic. If a recurrent homeomorphism on the torus, the annulus or the M$\ddot{o}$bius strip or the Klein bottle has a periodic point, then the homeomorphism is periodic. 
\end{satz}
Proof: See corollary [4.2] and remark [4.3] in [4].{\hfill $\Box$ \vspace{2mm}}

Proof of the theorem: 
Let $F:\stackrel{\circ }{ \Sigma } \to \stackrel{\circ }{ \Sigma }$ denote our return map.
If we show that our extension is recurrent, then by theorem [5.3] and our assumptions we know
that $F$ is periodic and every orbit of the flow is periodic, thus  we only need to show that $F$ is recurrent. \\ If the section $\Sigma$ is a  surface with boundary, then  $\stackrel{\circ }{ \Sigma }$ denotes the interior. 
Choose a finite number (say $j=1, \cdots , N$) of connected, compact, orientable surfaces $\Sigma_j$ that are diffeomorphic to discs and such that $\bigcup_j \Sigma_j= \Sigma $. Choose, for every
$j$, a sequence of  connected, compact, orientable surfaces $K_{j,i}$ such that
$K_{j,i} \subset \stackrel{\circ }{ K}_{j,i+1}$ and $\bigcup_i K_{j,i}  $ is an open set.
We construct $K_{j,i}$ such that given any compact set 
$C \subset \stackrel{\circ }{ \Sigma }$ we have, for a large $I$, that  $C \subset \bigcup_j K_{j,I}  $.
Moreover, we suppose  $\bigcup_{i,j} K_{j,i}  = \stackrel{\circ }{ \Sigma } $ and
$\bigcup_{j} K_{j,i}$ is connected. \\ For every $K_{j,i}$ we can find a tubular neighbourhood
 $U_{j,i}$ diffeomorphic $K_{j,i} \times (-1,1)$ via a diffeomorphism $\tau_{i,j}$. We  say  a point in $\tau_{i,j}^{-1} (K_{j,i} \times (-1,0))$ lies above $K_{j,i}$ and a point in 
$\tau_{i,j}^{-1} (K_{j,i} \times (0,1))$ lies below $K_{j,i}$. We set $U_{\ast, i,j}:=\tau_{i,j}(K_{j,i} \times (-1,0))$ and $U_{ i,j}^{\ast}:=\tau_{i,j}(K_{j,i} \times (0,1))$. We can suppose  that
$U_{\ast, i,j}\subset U_{\ast, i+1,j}$ 
and $U_{ i,j}^{\ast}\subset U_{ i+1,j}^{\ast}$. Since $\Phi$ is transversal to the section and without
singularities, we conclude that  $\bigcup_j K_{j,i}  $ is an orientable surface. Therefore we can define in tubular neighbourhood of $\bigcup_j K_{j,i}  $ what lies below and above. \\ Analogously we can define again  for $K_{j,i}$ and large $M$ the constant 
$$s(j ,i, M):=\inf\{|t_0-t_1|\,  |\,  -M<t_0< t_1 < M,  \Phi_{t_0}(v)\in \Sigma , 
\Phi_{t_1}(v)\in \Sigma , v\in K_{j,i} \}.$$ We can also define analogously 
the constant $\nu(j ,i, M, r) > 0 $. Given $ v\in \stackrel{\circ }{ \Sigma } $, then for every integer  $n$ we define $t(n,v)$ as the 
following:  \\ $t(n,v)$ is the unique element of $\mathbb{R}$ such that
 $\Phi_{t(n,v)}(v)=F^n(v)$ and 
$t(0,v)=0$. Moreover given $v\in \stackrel{\circ }{ \Sigma }$ and $t\ge 0$ we define:
$$P(v,t):= \max\{n \ge 0 \quad |\quad  t(n,v) \le t \}. $$ The counterpart
of lemma [3.4] will be the following:

\begin{lem}
Fix $v_j\in \bigcap_i K_{j,i}$.
There exist  sequences $T_i\to \infty$, $0< \zeta_i\to 0$ and $0 < \beta_i \to 0$
such that $\Phi_{T_i}(v_1) \in \stackrel{\circ }{ \Sigma }$  and for all $j\in \{1, \cdots N \}$ we have:
\begin{enumerate}
\item $\Phi_{T_i+ \zeta_i}(v) \in U_{ i+1,j}^*$ for all $v\in K_{j,i}$
\item $\pi(\Phi_{T_i- \zeta_i}(v)) \in U_{*, i+1,j}$ for all $v\in K_{j,i}$
 \item $d(\Phi_{T_i+s}(v),v)< \beta_i$ for all $v\in K_{j,i}$ and $|s|< 2\zeta_i$
\item For every $v\in K_{j,i}$ there exists an unique intersection point of $\stackrel{\circ }{ \Sigma }$
and $\Phi_{T_i+s}(v)$ where $s$ ranges over $|s|< 4\zeta_i$
\end{enumerate}
\end{lem}
Proof: Analogous to lemma [3.4].{\hfill $\Box$ \vspace{2mm}}

The proof is now easy. Set $p_{j,i}:=P(v_j, T_i+\zeta_i)$. The functions
$$ G_i: A\to \mathbb{N}$$ defined by $G_i(v)=P(v, T_i+\zeta_i)$ are constant on $K_{j,I_0}$ 
if $i\ge I_0$, but since $\bigcup_i K_{j,I_0}$ is connected, we know that $p_{j,i}$
is independent from $j$, and so $p_{j,i}=p_i$. We conclude $d(F^{p_i}(v),v)=d(\Phi_{t(v,p_i)}(v), v)<\beta_i$,
thus $F$ is paracompact-recurrent on  $\stackrel{\circ }{ \Sigma }$. Since the collapsed  boundary components  are fixed points, we 
conclude from lemma [2.6] that $F$ is recurrent on $X$.{\hfill $\Box$ \vspace{2mm}}

\begin{ko}
If $\Phi$ is a $C^{\infty}$ equicontinuous flow without singularities on a three dimensional manifold that admits a global surface of section $\Sigma $ and has at least three distinct periodic orbits, then the flow is pointwise periodic. 
\end{ko}

\section[]{Noncompact Manifolds}

In this section, we regard  noncompact manifolds, therefore two metrics that generate the same topology may be non equivalent. We can prove some fact about pointwise 
equicontinuous geodesic flows on noncompact surfaces if we restrict to a canonical metric for geodesic flows. In our case this should be the Sasaki metric.\\ Let $M$ denote a surface (maybe non compact).
$d$  will always denote the induced metric on $M$ of the Riemannian metric and $\tilde{d}$
the induced metric on $SM$ of the Sasaki metric. Note that we have by construction of the Sasaki metric $\tilde{d}(v,w)\ge d(\pi (v), \pi (w) )$.

\begin{defin}
 A system $(X,T)$ is called   pointwise equicontinuous  (pointwise regular) if for all $\epsilon > 0$ and all $x\in X$ there exists an $\delta (\epsilon , x) > 0$ such that for all $y$ with 
$d(x,y) < \delta (\epsilon , x) $ we have $d(xt,yt)< \epsilon$ for all $t \in T$.
\end{defin}

On compact metric spaces, pointwise equicontinuity implies equicontinuity, 
but on noncompact metric spaces this is not true and pointwise equicontinuity
is not independent from the metric that generates the topology.
Moreover, given any compact set $K$ of $X$,  we can find an $\delta (\epsilon , x) > 0$
in the definition above  that is independent from $x\in K$.

\begin{lem}
Let $(\Phi , SM, \tilde{d} )$ be a geodesic flow that is pointwise equicontinuous.
For all compact  sets $K$ there exists a number $C(K)$ such that for all $v,w \in K$ and $t\in \mathbb{R}$ we have $d(\gamma_{v_0}(t), \gamma_w(t)) \le C(K)$. 
\end{lem}

Proof:  Choose first an open and bounded set $O$ that contains $SM|K$. Choose for
 $\overline{O}$
a constant  $\delta (\epsilon) > 0$ such that if for $v,w\in SM|K$ we have $d(v,w)< \delta (\epsilon)$ 
then $d(\Phi_t (v), \Phi_t (w)) < \epsilon $. Cover $SM|K$ with $N=N(\epsilon , K)$
balls of radius smaller than $\delta (\epsilon)$ such that the union of these balls is connected
and lies in $O$. We conclude 
$$ d(\gamma_v (t), \gamma_w (t)) \le N\epsilon := C(K)$$
for all $v,w \in SM|K$ and $t\in \mathbb{R}$.{\hfill $\Box$ \vspace{2mm}}

\begin{hl}
If $M$ is compact and the geodesic flow $(\Phi , SM)$ is equicontinuous, then $\pi_1(M)$ is finite.
\end{hl}
Proof: $(\Phi , SM)$ is equicontinuous with respect to the Sasaki metric and therefore
one can conclude that the lifted geodesic flow on the universal covering $\tilde{M}$ is equicontinuous with respect to the lifted Sasaki metric 
(denoted with $\tilde{d}$). Since $M$ is compact, we conclude that for all $r>0$ there exists a number $Z(r)$ such that for any point  $x$  of $M$ the sphere $S_xM$ can be covered by $Z(r)$ balls of radius $r$ (with respect to the metric $d$).  Given any point $x$ and assume $M$ is not compact. Choose $v\in S_xM$ such that $\gamma_v : \mathbb{R}^+ \to M$ is a ray and a sequence $t_i\to \infty$. Set 
$v_i=-\dot{\gamma_{v_i}}(t_i)$. We have
by the proof of the lemma above and the existence of $Z(r)$ for an $\epsilon >0$
$ d(\gamma_{v_i}(t_i), \gamma_{-v_i}(t_i)=\gamma_{v}(2t_i)) \le Z(\delta (\epsilon) ) \epsilon$, 
hence $2t_i=d(x=\gamma_{v_i}(t_i), \gamma_{v}(2t_i) ) \le Z(\delta (\epsilon) ) \epsilon $,
but $t_i$ grows.{\hfill $\Box$ \vspace{2mm}}

\begin{ko}
If $M$ is noncompact and the geodesic flow  $(\Phi , SM)$ is pointwise equicontinuous with respect to the Sasaki metric, then the following holds:
\begin{enumerate}
\item There is no minimal geodesic (also called line) in $M$.
\item $Diam (\partial B(x, r)): = \sup \{ d(x,y) | x,y\in \partial B(x, r) \} $ is bounded by a constant $C(x)$ for all $x$. The constant $C(x)$ can be chosen uniformly on compact sets.
\end{enumerate}

\end{ko}
Proof:  If $\gamma$ is minimal, then for some constant $C$ 
$$2t = d(\gamma (t), \gamma (-t)) \le C, $$
 and therefore $\gamma$ is not minimal. \\ If $Diam (\partial B(x, r))$ is not bounded,
choose sequences $a_i$ and $b_i$ such that $a_i, b_i \in \partial B(x, r_i) $ and $d(a_i, b_i)\to \infty$ where $r_i \to \infty$. Choose for $a_i$ a sequence $v_i\in S_xM$ such that 
$\gamma_{v_i}:[0,d(x, a_i)]$ is a minimal geodesic segement that starts in $x$ and ends in $a_i$. W.l.o.g. we suppose that $v_i\to v$ and, therefore $\gamma_v$ is a ray and from equicontinuity we conclude that $d( \gamma_v (r_i) , a_i) \to 0$. We repeat the construction for $b_i$
to get a ray $\gamma_w$ such that $d( \gamma_w (r_i) , b_i) \to 0$. The flow is equicontinuous
and therefore we have -by lemma [6.2]- $d(\gamma_v (r_i), \gamma_w (r_i)) $ bounded 
and therefore $d(a_i, b_i)$ is bounded.  It follows from the proof that $C(x)$ can be chosen uniformly on
compact  sets.{\hfill $\Box$ \vspace{2mm}}

Of course, one can conjecture that theorem [4.4] holds in higher dimensional cases, but it seems that there are no tools to prove this conjecture. Corollary [6.3] holds for P-manifolds (see theorem 7.37 in [3]) and  using the Morse-Index theorem, one can 
see that noncompact manifolds with strictly positive sectional curvature have no line,
thus there are some reasons to conjecture this. \\  We now discuss if there exists an  equicontinuous (or even a pointwise equicontinuous) geodesic flow with respect to the 
Sasaki metric on a noncompact surface. \\ Here is a subresult:

\begin{hl}
Let $(M, g)$ be a  Riemannian manifold of dimension 2 and suppose that the geodesic flow
$(\Phi , SM)$ is pointwise equicontinuous with respect to the Sasaki metric, then $M$ is homeomorphic to the plane. 
\end{hl}

Proof: Our proof is based on the following nice theorem 
\begin{satz}
Every surface is homeomorphic to a surface formed from a sphere $S$ by first removing a closed totally disconnected set $X$ from $S$, then removing the interiors of a finite or infinite sequence $D_i$ of non overlapping closed discs in $S-X$ and finally suitable identifying the boundaries of these discs in pairs. It may be necessary to identify the boundary of one disc with itself to produce a "cross cap". The sequence $D_i$ approaches $X$ in the sense that, for any open set $U$ in $S$ containing $X$, all but a finite number of  the $D_i$ are contained in $U$.
\end{satz}
Proof: See [7].{\hfill $\Box$ \vspace{2mm}}

 $d$ is the metric induced from the Riemannian metric. 
We first prove that $M$ has  genus zero! Assume that $M$ is not diffeomorphic to $S^2-X$.
Choose a curve $\beta$ that is not contractible and lies on a "handle" or a "crosscap" such that we can find a compact set $O$ that contains $\beta$ and is diffeomorphic to a closed  disc where a cylinder or a crosscap is glued in. Let $[\beta]$ denote the free homotopy class. If there is a sequence $\beta_i$ such that $\beta_i \in [\beta]$ and $L(\beta_i) \to 0$,
then, since $\beta_i$ always intersects $O$, we conclude that $\beta_i$ is contradictible
for large $i$. Now let $\beta_i$ be  a sequence such that $\beta_i \in [\beta]$ and $L(\beta_i) \to \inf_{c \in [\beta]} L(c) > 0$. Again we know that $\beta_i$ always intersects $O$.
If $\beta_i$ lies in a compact subset  of $M$, then we know that a subsequence converges to a
closed geodesic, but it is clear that $\beta_i$ lies in a compact subset $K$ of $M$, since 
otherwise
there is a point $x_i\in \beta_i$ that tends to infinity, hence $L(\beta_i)$
is not bounded. Hence we have found a closed geodesic $\beta$. Since there exists a ray starting at a point of $\beta$, we conclude from lemma [6.2] that the ray does not tend to infinity, thus we get a contradiction. \\ We prove now that $X$ is just a simple point.
We endow $S^2$ with a metric $d_0$ that generates the standard topology of $S^2$.
Let $-\infty$ and $\infty$ two different points of $X$. Choose  sequences $a_i$ and $b_i$ such that $a_i\to -\infty$ and $b_i\to \infty$ with respect to $d_0$. We want to show that $d(a_i, b_i)$  tends to infinity for a subsequence. Suppose that $d(a_i, b_i)$ is bounded. Let $\delta_i: [0, d(a_i, b_i)] \to M$ be a minimal geodesic segment from $a_i$ to $b_i$. Let $K_i$ denote the image of $\delta_i$. Define $$\lim (K_i):=\{y\in S^2\quad | \quad  \exists x_i\in K_i : d_0(x_i,y) \to 0 \}.$$ It is easy to see that $K:=\lim (K_i)$ is closed
and connected.\\ Assume a point of $M$ lies in $K$. Then  w.l.o.g. a tangent vector $v_i$ of $\delta_i$ converges to a vector $v$ of $SM$. Since we suppose that $d(a_i, b_i)$ is bounded by a constant $C$, then $\gamma_v$ will be a geodesic that meets $-\infty$ and $\infty$
on the intervall $[-2C, 2C]$, since $\delta_i[0,C] $ converges to a segement of $ \gamma_v[-2C,2C] $ with respect to $d$. Hence $K$ is a subset of $X$, but this means $K$ is reduced to a point, thus $d_0(a_i, b_i) \to 0$ and therefore we have a contradiction to the fact that $d(a_i, b_i)$ is bounded. Thus given any sequence $a_i \to -\infty$  and $b_i \to \infty$ then $d_1(a_i, b_i)$ tends  to infinity for a subsequence. \\ We now construct  sequences $n_i \to -\infty$  and $m_i \to \infty$ such that $d(n_i, m_i)$ is bounded and therefore we derive a contradiction
to the fact that $X$ contains more than one point. Define 
$$\omega (v):= \{   y\in S^2\quad | \quad  \exists t_i\to \infty :
d_0(\gamma_v(t_i),y) \to 0  \}. $$
It is easy to see that $\omega (v)$ is closed and connected in $S^2$. If, for a ray
$\gamma_v$, the set  $\omega (v)$ is  not subset of $X$, then there is a sequence 
$t_i \to \infty $ such that $\dot{\gamma_v}(t_i)$ converges to a vector $v^*$ of $SM$
and hence $\gamma_{v^*}$ will be a minimal geodesic, since $\gamma_v$ is a ray. This contradicts corollary [6.4]
and therefore $\omega (v)$ is a subset of $X$, hence for a ray $\gamma_v$, the set 
$\omega (v)$ is reduced to a point of $X$ ($X$ is totally disconnected).\\  Take  sequences $a_i \to -\infty$  and $b_i \to \infty$. For choose $a_i$ a sequence $v_i\in S_xM$ such that 
$\gamma_{v_i}:[0,d(x, a_i)]$ is a minimal geodesic segement that starts in $x$ and ends in $a_i$. $v_i$ tends to a vector $v$. The curve $\gamma_v$ will be a ray and therefore
$\omega (v)$ will be a point. We show that $-\infty = \omega (v)$. Note that $d(a_i,\gamma_v (t_i)) \to 0 $, since our flow is equicontinuous in $x$. Suppose that $d_0(\gamma_v (t_i), q)\to 0$ for a point $q\in S^2$. If $q\not= -\infty$, then choose a minimal geodesic segment
 $\delta_i : [0, d(a_i,\gamma_v (t_i))] \to M$ that starts in $a_i$ and ends in $\gamma_v (t_i)$.  We denote the image of $\delta_i $ with $K_i$. Again we conclude as above, that $\lim (K_i)$  will be a connected subset of $X$, and therefore $q= -\infty$. \\ Therefore we get
a ray $\gamma_v$ that starts in $x$ and converges to $-\infty$. Repeating this
constructions will generate a  ray $\gamma_w$ that starts in $x$ and converges to $\infty$.
From lemma [6.2] we know that $d( n_i:=\gamma_v(i), m_i:=\gamma_w(i))$
is bounded.{\hfill $\Box$ \vspace{2mm}}

We end this paper by showing that at least a pointwise equicontinuous geodesic flow 
with respect to a metric  $d_1$  exists. The metric 
$d_1$ is not equivalent to the induced metric of the Sasaki metric of $h$ ($h$ will be defined later).  Let $g_0$ be the standard metric on $\mathbb{R}^2$, $d_0$ the metric induced from $g_0$ and identify $ \mathbb{R}^2$ with $ \mathbb{C}$.    Choose a diffeomorphism $f:[0,\infty) \to [0,\infty)$
such that $f|[0, \frac{1}{2})= Id$ and $f(t)=\exp (t)$ for $r \ge 1$. 
Let $F$ denote a diffeomorphism 
from $\mathbb{R}^{> 0} \times S^1$ to  $\mathbb{R}^2 -\{ 0 \}$
defined by $F(r, \phi) =(f(r),\phi) $ where $(r, \phi)$ are the polar coordinates. The pull-back of the metric $g_0$ under $F$ defines a new metric $h$ on $\mathbb{R}^2$ and its geodesic flow is complete. Regard the metric $d_1$ on $S\mathbb{R}^2$ defined by
$d_1((x,v),(y,w))= \| x-y \| + \| v-w\| $. We show that the geodesic flow of 
$(\mathbb{R}^2, h)$ will be pointwise equicontinuous with respect to $d_1$.\\ Consider the coordinante system $F: \mathbb{R}^2-\{0\} \to \mathbb{R}^2-\{0\} $. Given a vector 
$v$ at a point $x$, we regard the geodesic $g_{x,v}(t)=x+tv $ of the metric $g_0$. Let us fix 
a point $x=a_0+ib_0$ and a vector $v= a_1+ib_1$. Since our metric is  invariant with respect
to revolutions, we can assume that $a_1\not= 0$.  For a small $\epsilon > 0$ choose a $M> 0$ and a $\delta (\epsilon )> 0 $ such that  
$d_0(\dot{g}_{x,v}(t), \dot{h}_{y,w}(t) )< \epsilon $ for $t\in [-M,M]$ and $d((x,v),(y,w))< \delta (\epsilon )$. If $M$ is large and $\delta (\epsilon )> 0 $ small enough, then  the geodesic $h_{y,w}(t)$ lies outside the compact set 
$F([0,2], \mathbb{R} )$ for $t\notin [-M, M] $. \\ Note that $|g(t)|= g(t) \cos \phi(t) +ig(t)\sin \phi(t) $
 ,where \\ $\phi(t) = \arctan (\frac{b_0+tb_1}{a_0+ta_1})$ and $|g(t)|^2= (a_0+ta_1)^2+(b_0+tb_1)^2$ . Note that 
$$F^{-1}= g^*_{x,v}(t):=  (f^{-1} (|g_{x,v}(t)|), \arctan (\frac{b_0+tb_1}{a_0+ta_1}) )  $$
will be a geodesic of our manifold $(\mathbb{R}^2, h)$. \\ We show that if $\epsilon > 0$ is small enough, the distance $d_1(\dot{g}^*_{x,v}(t), \dot{h}^*_{y,w}(t) ) $  remains arbitrary small for all $t$. A computation shows that for $t\notin [-M, M] $ we have 
$$\dot{g}^*_{x,v}(t)= (\frac{a_1(a_0+ta_1)+ b_1(b_0+tb_1)}{|g(t)|^2}, 
\frac{b_1(a_0+ta_1)+ a_1(b_0+tb_1)}{(1 + ( \frac{b_0+tb_1}{a_0+ta_1} )^2 )
(a_0+ta_1)^2} ).$$
Hence for a small $\epsilon >0$ and large $M$, we have $\|\dot{g}^*_{x,v}(t) - \dot{h}^*_{y,w}(t) \|$  small for $t\notin [-M,M]$.
Let $h_{y,w}(t)=y+tw= n_0+im_0+t(n_1+im_1)$, then one can compute
that $$\|h_{y,w}(t)- g_{x,y}(t) \|\le \|f^{-1}(|g(t)|) -f^{-1}(|h(t)|)  \| $$
 $$+ \| \arctan (\frac{b_0+tb_1}{a_0+ta_1} ) - \arctan (\frac{m_0+tm_1}{n_0+tn_1})    \| .$$
In the same way we  conclude that for  a small $\epsilon > 0$ and large $M$, we will have the second term   small for $t\notin [-M,M]$.  A computation shows that 
$$ \|f^{-1}(|g(t)|) -f^{-1}(|h(t)|)  \| = 
\|\frac{1}{2} \ln (\frac{ \frac{a_0^2+b_0^2}{t^2}+ \frac{2(a_0a_1+b_0b_1)}{t}+a_1^2+b_1^2 }{ \frac{n_0^2+m_0^2}{t^2}+ \frac{2(n_0n_1+m_0m_1)}{t}+n_1^2+m_1^2 } ) \| ,  $$
but this term will be small for  a small $\epsilon > 0$, large $M$ and $t\notin [-M,M]$. Therefore we know that the geodesic flow is equicontinuous with respect  to the metric $d_1$.
The metric is not equivalent to the induced metric of $h$, since the Riemannian distance 
between two geodesics of different directions  grows.

\section[]{Acknowledgements}

The author thanks Gerhard Knieper for asking the author the nice question if strong
recurrence behavior of the geodesic flow  implies that all geodesics are closed.
This inspired  the author to write this article. The author also thanks
Henrik Koehler, Yoro Olivier, Serge Corteyn for the corrections
and remarks.

\end{document}